\documentclass[9pt,letterpaper]{IEEEtran}
\usepackage{amsmath}
\usepackage{amssymb}  
\usepackage{palatino, url, multicol}
\usepackage[numbered]{mcode} 
\usepackage{graphicx}

\usepackage{url} 

\begin{document}

\title{Effect of multiple degrees of ambivalence on the Naming Game}
\author{
    \IEEEauthorblockN{Yosef Treitman, Chjan Lim}\\ 
    \today
}

\maketitle

\begin{abstract}

We examine a modified Naming Game in the mean field where there are
multiple degrees of ambivalence.  Once an agent in one state fears
an opinion one way or another, he or she moves one step in the
appropriate direction. In the absence of zealots, the two consensus
states are stable steady states and the uniform distribution is an
unstable steady state.  With zealots for one opinion only, there is
a critical value below which there are three steady states and above
which there is only one.  Consensus in favor of the zealots' opinion
is the steady state that always exists, and is stable.  The second
steady state is the uniform distribution in the absence of zealots,
and moves away from the zealots' opinion as the number of zealots
increases.  This state is unstable.  The last steady state starts at
consensus against the zealots, and moves toward the zealots' opinion
as the number of zealots increases.  This state is stable.  When
zealots are added on both sides, the "beak" pattern observed for the
Naming Game remains, with the region of multiple steady states
growing with the addition of more intermediate states.

\end{abstract}



\section{Introduction}

The science of social interaction is a very alluring field.  To
understand how humans influence each other provides the opportunity
to control how public opinion shapes itself.  The difficulty is that
the human brain is too complex to allow for a simple and clearly
accurate model of human opinion formation.  As such, a number of
approximations have been proposed, each with its own shortcomings.
Of interest is the rather simple voter model.  Here, each person, or
agent, has one of two conflicting opinions, A or B.  As the agents
talk to each other, they try to convince each other of their own
opinions.  One agent is chosen as the speaker and the other is
chosen as the listener.  If the speaker and listener agree, then
there is no change in the opinion state, but if they initially
disagree, then the listener is converted to the speaker's opinion.
[2, 4, 5, 9, 11, 12, 17, 18] Of particular interest is the parameter
$m$, the magnetization of the system.  The magnetization is the
expected poll result and is calculated by $\frac {\rho_A} {\rho_A +
\rho_B}$, where $\rho_A$ represents the density of people in favor
of opinion A and $\rho_B$ represents the density of people in favor
of opinion B.  Note that, for the voter model, $\rho_A + \rho_B =
1$.  In this case, the variables $\rho_A$ and $\rho_B$ are
martingales, and given enough time, the system will almost surely go
to one of two absorbing states, that of full consensus at A and that
of full consensus at B. [2, 4, 5, 9, 11, 12, 18]
\newline

A more interesting model is the Naming Game model, or Binary
Agreement Model. Here, the possibility of ambivalence is accounted
for. In addition to the two extremist states A and B, there is the
state AB, which represents ambivalence.  In this model, agents in AB
are no more or less likely to speak or be spoken to than anyone
else. If the chosen speaker is in the ambivalent state, then he or
she will subconsciously choose one of the two opinions.  If the
listener is in the ambivalent state, then he or she will change to
whichever opinion he or she hears. [1, 3, 4, 5, 6, 7, 8, 13, 14, 15,
16, 17, 18, 19] While there are various psychological studies on the
effect of voicing an opinion on an ambivalent speaker, this is
ignored in the Listener-Only Naming Game model. [4, 18] In this
model, the values $\rho_A$, $\rho_B$, and $\rho_{AB}$ are no longer
martingales.  In fact, once one side begins to dominate, (this
dominance can be determined by the magnetization, which is now
calculated by $m = \frac {\rho_{AB} + 2\rho_A}{2}$) this dominance
will tend to grow until consensus is reached.  The consensus points
are stable steady states while perfect balance between $\rho_A$,
$\rho_B$, and $\rho_{AB}$ is an unstable steady state. [1, 3, 4, 5,
6, 7, 8, 13, 14, 15, 16, 17, 18, 19]
\newline

A further generalization allows for the agents to be leaning in one
direction or another to varying degrees without full commitment to
one idea or the other.  If a listener hears one opinion or the
other, he or she moves one step in the direction of the opinion he
or she hears. [18] The parameter $K$ represents the number of times
that an agent convinced of opinion B needs to hear opinion A in
order to be convinced of opinion A, and vice versa.  The population
densities in each state are now denoted $\rho_0, \rho_1, ...\rho_K$,
where $\rho_0$ represents commitment to opinion B, $\rho_K$
represents commitment to opinion A, and the other states are the
various ambivalent states.  The states themselves (as opposed to the
populations of the states,) are denoted $N_0, N_1... N_K$.  An agent
in state $N_i$ will speak in favor of opinion A with probability
$\frac {i} {K}$ and opinion B with probability $\frac {K-i}{K}$.
Note that the Voter Model corresponds to $K = 1$ and the Naming Game
model corresponds to $K = 2$. [18] \newline

One other factor involved is the presence of zealots in favor of one
opinion of the other.  A zealot operates with a motive that ignores
all logic and thus can never be convinced to change his or her
opinion. No amount of convincing can turn a non-zealot, or normal
agent, into a zealot.  The presence of these zealots affects the
dynamics of the system as well as the possible long-term outcomes.
[1, 3, 4, 9, 11, 12, 13, 14, 17, 18]

\section{The No-Zealot Case}

The behavior of a network obeying one of these models can be
characterized by its long term behavior in the mean field.  That is
to say, what possible steady-states are there and which of them are
stable. [20] It has already been shown that for the standard Naming
Game model, the two consensus points are stable steady states and
that the uniform distribution between $\rho_0$, $\rho_1$, and
$\rho_2$ is an unstable steady state. [1, 3, 4, 5, 6, 7, 8, 13, 14,
15, 16, 17, 18, 19] As it turns out, the analogous result is true
for any value of $K \geq 2$. Namely, the state $\rho_0 = 1$, $\rho_i
= 0$ for $i \neq 0$ is a stable steady state, as is the state
$\rho_K = 1$, $\rho_i = 0$ for $i \neq K$. The only other steady
state is at $\rho_i = \frac {1}{K+1}$ for all i, and this state is
unstable.
\newline

The outline of the proof is as follows:  It has been previously
shown that a necessary condition for a steady state is that $\frac
{\rho_i}{\rho_{i-1}}$ is constant (if finite,) for $i = 1, 2... K$.
Furthermore, this quotient is equal to $\frac {m}{1-m}$, where $m$
is the magnetization, or the expected poll result.  The
magnetization is calculated with the formula $m = \sum_{i = 1}^{K}
(\frac {i}{K} \times \rho_i)$.  It can be shown that this is a
necessary and sufficient condition for a distribution to be a steady
state.  Thus, selecting a value of $m$ between 0 and 1 fully
determines the distribution by fixing the ratios of population
densities of different states.  However, the distribution has a
magnetization of its own, and this may or may not be the same as the
magnetization used to generate the distribution.  If it is the same,
then it is a steady state, otherwise it is not.  It can be shown
that for $m = 0$, $m = 0.5$, and $m = 1$, the resulting distribution
has the same magnetization as the one used to generate it.  For
$0<m<0.5$, the resulting magnetization is less than the value of $m$
used to generate the distribution, and for $0.5<m<1$ the resulting
magnetization is greater than the one used to generate the
distribution.\newline

To address stability, we first look at the states corresponding to
$m = 0$ and $m = 1$.  Because the setup exhibits symmetry, we only
need to show stability for one point, and the other point should
have the same level of stability.  It should be noted that the only
distribution with $m = 0$ is $\rho_0 = 1$, $\rho_i = 0$ for $i \neq
0$.  For distributions near this one, $m$ tends to decrease over
time.  Thus, the distributions tend toward the minimum value of $m$,
which is $m=0$, and there is only one distribution with that
magnetization value.  The point at $m = 0.5$ has a different
property.  Assuming a geometric distribution with $m$ slightly
greater than 0.5, the value of $m$ will tend to increase, and if $m$
is slightly less than 0.5, $m$ will tend to decrease.  Thus, if the
distribution is slightly perturbed from the steady state, the
perturbation will grow larger and larger, resulting in instability.
A more detailed proof appears in the appendix.

\section{Steady States for the Unilateral Zealot Case}

We now examine the case where, in addition to the normal agents,
there are agents with unshakable support for one of the two
opinions.  We perform calculations assuming the zealots support
opinion B, and use symmetry to infer the corresponding results for
zealots in favor of opinion A.  We denote the population density of
zealots in favor of opinion B as $\rho_B$.  It is worth mentioning
that the notion of long-term behavior takes on multiple meanings in
this case.  If there are no zealots in favor of opinion A, but there
is at least one in favor of opinion B, then there is only one
absorbing state, namely, that of full consensus at opinion B.
Furthermore, it can be shown that this consensus state will be
reached with probability 1. [1, 3, 4, 13, 14, 17, 18, 19, 20, 21,
22] So in a sense, the only truly stable steady state is that of
consensus at opinion B. However, it is realistically possible for
the system to hover in the neighborhood of a particular distribution
for a relatively long period of time before drifting to the
consensus state. [21, 22] It is these particular distributions which
we are looking for.  We will call these opinion states "metastable"
states. Additionally, we look for the conditions under which these
metastable states can exist.\newline

In the case of the Voter Model, it is rather simple.  The
magnetization is a supermartingale, with martingality holding only
for the consensus state. [3, 4, 9, 13, 17, 18] In the Naming Game
model, things get more interesting.  Consensus at B is always a
stable steady state. It has been shown that in the case of $\rho_B =
0$, consensus at A is also a stable steady state, and that this
stable state drifts toward lower values of $m$ as $\rho_B$
increases.  This is also true for the case K = 3.  Additionally,
when $\rho_B = 0$, the uniform distribution is an unstable steady
state.  This unstable state moves to one of higher magnetization as
$\rho_B$ increases. [3, 4, 13, 14, 17, 18, 19, 21] This has also
been found to be true in the case K = 3. To understand why this is,
one needs to understand the relationship between the magnetization
$m$ and the magnetization of the normal agents, $m_{normal}$.  We
calculate the magnetization of the normal agents by the formula
$m_{normal} = \frac {\sum_{i = 1}^{K} i \times rho_i}{K \times
(1-\rho_B)}$.  If we fix a value of $m$ and force the normal agents
into a geometric distribution obeying $\frac {\rho_i}{\rho_{i-1}} =
\frac {m}{1-m}$, then the value of $m_{normal}$ is determined.  It
can be shown that the number of zealots committed to opinion B can
be determined from the formula $\rho_B = 1 - \frac {m} {m_{normal}}$
As it turns out, the results regarding the existence and
metastability of steady states are true as long as $\frac {d^2
\rho_B}{dm^2}$ is of one sign over the interval $(0.5,1)$ and
another over the interval $(0,0.5)$. Figure 1 shows that this is the
case for K = 2, 3, 4, and 10. Figures 2-5 examine the stability of
steady states in the cases K = 3 and K = 10.  The outline of the
proof is as follows:\newline

\begin{figure}
\centering
\includegraphics[width=3.0in,bb=2.0in 3.2in 6.5in 7.5in]{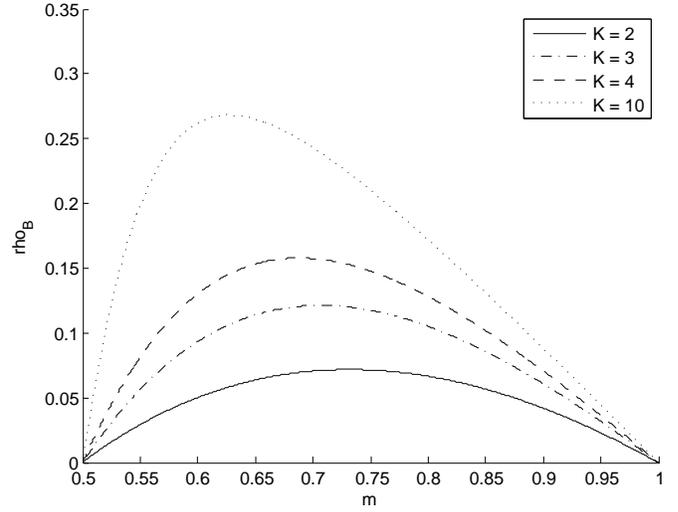}
\caption{The values of $\rho_B$ that will make a given value of m a
steady state for K = 2, 3, 4, and 10.  Note that, in each case, the
value of $\rho_B$ increases up to a point, then decreases.  Thus,
for each value of $\rho_B$ below a critical value, there will be
exactly 2 steady states in the interval [0.5, 1] and above this
critical value, there will not be any steady states in that
interval.}
\end{figure}

\begin{figure}
\centering
\includegraphics[width=3.0in,bb=2.0in 3.2in 6.5in 7.5in]{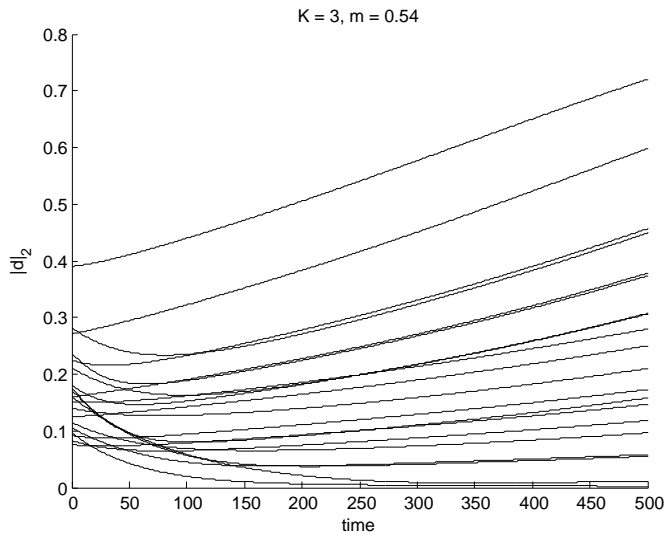}
\caption{Setting K = 3, we pick initial conditions with random small
perturbations off the intermediate steady state, and track the
2-norm of d, the difference between the perturbed state and the true
steady state.  Over the course of time, the difference tends to
grow, showing that the steady state is unstable.}
\end{figure}

\begin{figure}
\centering
\includegraphics[width=3.0in,bb=2.0in 3.2in 6.5in 7.5in]{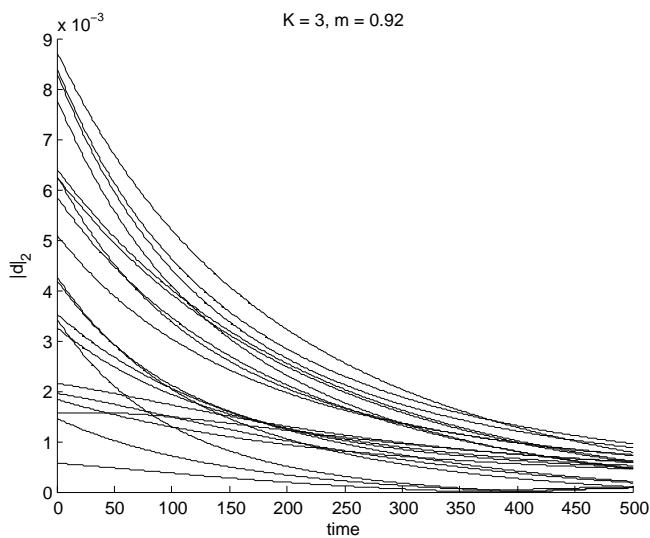}
\caption{Similar to figure 2, this time we take small perturbations
off the equilibrium state of highest magnetization.  The difference
always decreases to 0, showing this equilibrium state to be stable.}
\end{figure}

\begin{figure}
\centering
\includegraphics[width=3.0in,bb=2.0in 3.2in 6.5in 7.5in]{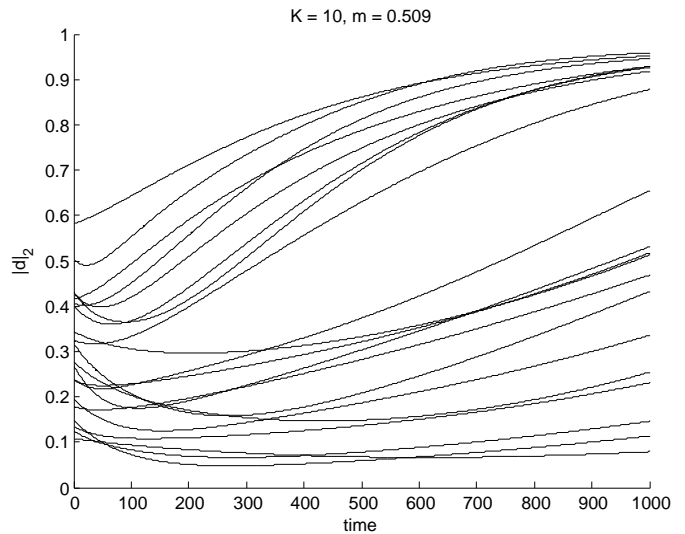}
\caption{Similar to figure 2, but with K = 10 and the intermediate
steady state. Because the difference grows over time, the steady
state is unstable.}
\end{figure}

\begin{figure}
\centering
\includegraphics[width=3.0in,bb=2.0in 3.2in 6.5in 7.5in]{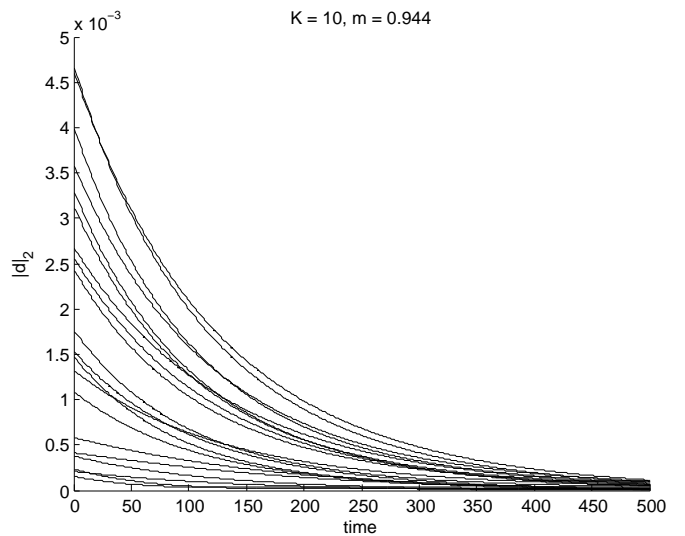}
\caption{Similar to figure 3, but with K = 10 and the
high-magnetization steady state. Because the difference shrinks over
time, the steady state is stable.}
\end{figure}

Similar to the no-zealot case, a necessary condition for a steady
state is the relationship $\frac {\rho_i}{\rho_{i-1}} = \frac
{m}{1-m}$.  A subtle difference between this case and the no-zealot
case is that $\rho_B$ is now a parameter, and while it does not
directly affect any of the $\rho_i$ values, it does directly affect
$m$.  In particular, as $\rho_B$ increases, $m$ decreases.  Thus, in
cases where $m_{normal} > m$, there is a suitable $\rho_B$ that will
enable the equation $\frac {\rho_i}{\rho_{i-1}} = \frac {m}{1-m}$ to
hold.  Recall the result that if a geometric distribution is
generated from a magnetization $m$ then the resulting magnetization
will exceed the generating magnetization for $0.5<m<1$.  Thus, for
$0.5<m<1$, there is a positive value of $\rho_B$ that will make the
distribution steady.  By symmetry, for $0<m<0.5$, there is a value
of $\rho_A$ that will make the distribution steady if $\rho_B = 0$.
It can be shown that this special value of $\rho_B$ is unique and is
a continuous function of $m$.  This means that, as long as the
special value of $\rho_B$ is of one concavity as $m$ varies, (in
this case, concave down,) the steady states will each move in one
direction as $m$ increases.  Because $\rho_B$ is bounded by 0 and 1,
there must be some maximum value of $\rho_B$ that yields a steady
state other than consensus at B.  The continuity of $\rho_B$ shows
that multiple steady states must meet at this maximum.  Empirical
evidence suggests that the stability of these steady states does not
change until they meet at this critical value of $\rho_B$. A more
detailed proof appears in the appendix.

\section{Tipping Points for the Case of Bilateral Zealots}

As we have already mentioned, stability, and thus, long-term
sustainability of non-consensus states depends on the number of
zealots.  This effect is even further pronounced when there are not
only zealots in favor of opinion B, but also in favor of opinion A.
[10, 11, 12, 17, 19] Examining the data in figures 6-9, which
examine the conditions under which multiple steady states exist,
several phenomena are noticed. First, if there are sufficiently few
committed agents total, there will be multiple steady states.
Second, if there are significantly more zealots in favor of one
opinion than the other, there is only a single steady-state.
Finally, if there are enough zealots, there will only be one steady
state. [10, 11, 12, 17, 19] This holds even when the zealots are
perfectly evenly divided between opinions A and B. [11, 17, 19] This
case is of particular interest, as it implies that, given enough
zealots, the equilibrium uniform distribution becomes a stable
steady state. [10, 11, 12, 17, 19] More interestingly, it becomes
the only stable steady state, and thus, given enough time, the
system should always be near this state.

\begin{figure}
\centering
\includegraphics[width=3.0in,bb=2.0in 3.2in 6.5in 7.5in]{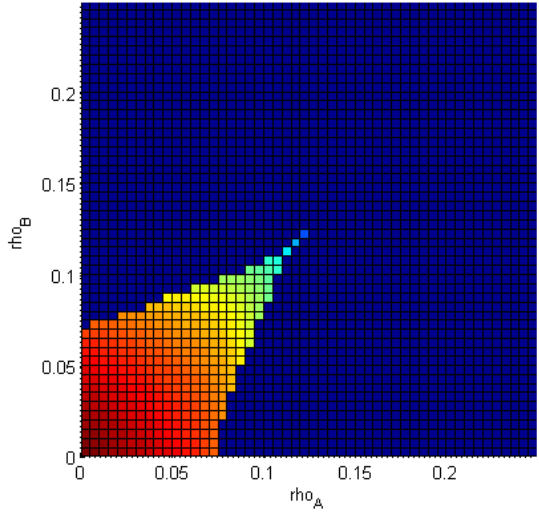}
\caption{An analysis of which combinations of $\rho_A$ and $\rho_B$
result in multiple steady states for K = 2.  The dark region is
where there is only one steady state.  Note the "beak" shape
consistent with [17].  Note that when $\rho_A$ = $\rho_B$, the dark
region begins at (0.125, 0.125).}
\end{figure}

\begin{figure}
\centering
\includegraphics[width=3.0in,bb=2.0in 3.2in 6.5in 7.5in]{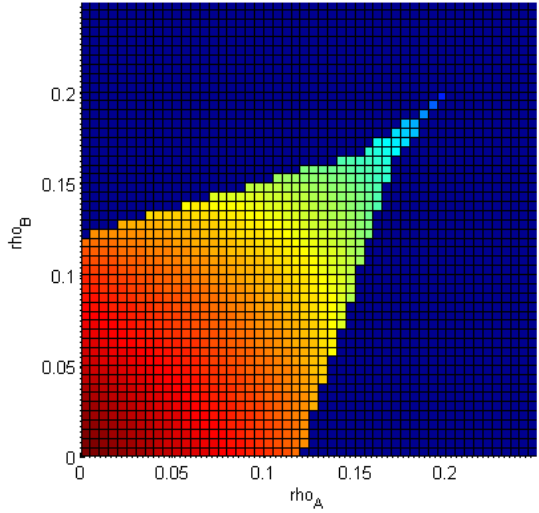}
\caption{Similar to figure 6, this graph analyzes the case of K = 3.
 The "beak" shape remains, and the cusp of the dark region is at (0.20, 0.20).}
\end{figure}

\begin{figure}
\centering
\includegraphics[width=3.0in,bb=2.0in 3.2in 6.5in 7.5in]{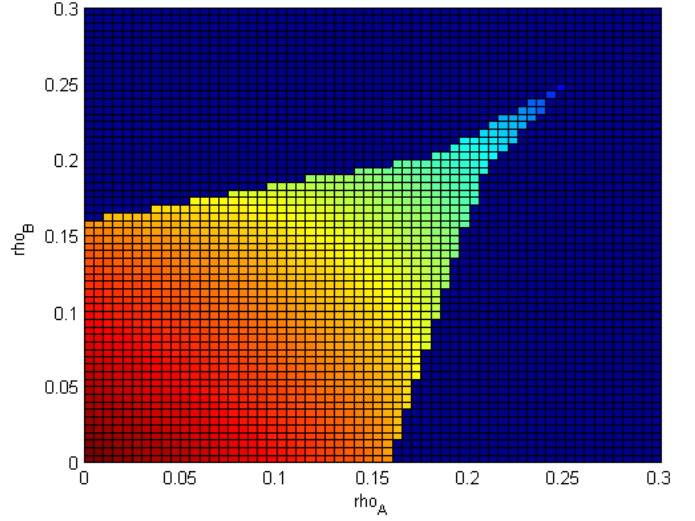}
\caption{This graph analyzes the case K = 4.  The cusp of the dark
region is now at (0.25, 0.25).}
\end{figure}

We find the critical value of $\rho_A$ and $\rho_B$ above which the
uniform distribution of normal agents is stable.  Let $\alpha$
represent the fraction of agents who are not zealots.  We also
assume that the system is near the uniform distribution.  Because
there is little drift here, we can assume that the system is near a
geometric distribution, as the value of $m$ would have been
approximately constant for a long time.  Because the system is
nearly uniform, we assume $r = 1 + \epsilon$.  This gives us the
following equations:

\begin{figure}
\centering
\includegraphics[width=3.0in,bb=1.5in 3.2in 6.5in 8.0in]{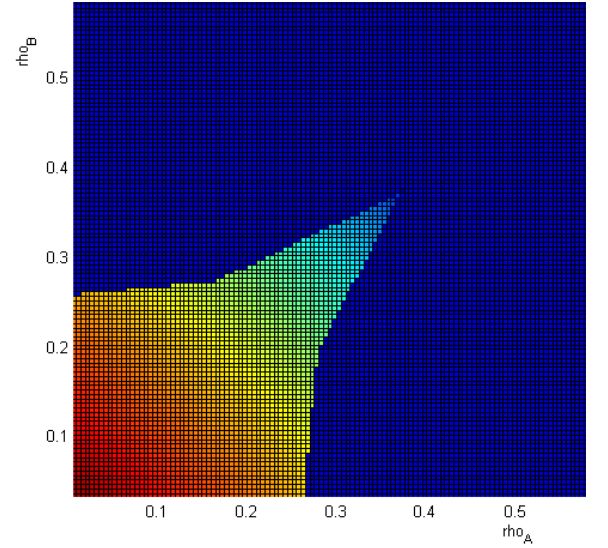}
\caption{This graph analyzes the case K = 10.  The "beak" pattern
still remains, and the cusp is at (0.375, 0.375).}
\end{figure}
\begin{align}
\rho_i = \rho_0(1 + \epsilon)^i\\
\end{align}
for $i = 0, 1, 2...K$

Note that, because $\epsilon$ is small, $(1 + \epsilon)^i \approx 1
+ i\epsilon$, giving us

\begin{align}
\rho_i = \rho_0(1 + i\epsilon)\\
\end{align}
for $i = 0, 1, 2...K$

To calculate $m$ for this distribution, we can first calculate $m$
for the zealots, then calculate $m$ for the normal agents, and take
a weighted average.  Because the zealots are evenly distributed
between $\rho_A$ and $\rho_B$, we know that, for this group of
zealots, $m_{zealots} = 0.5$.  For the normal agents, the
calculation is more involved.

\begin{align}
m_{normal} = \frac {\sum_{i=0}^{K} \frac{i+i^2 \epsilon}{K}} {K+1+
\frac {K^2 + K}{2}\epsilon}
\end{align}

Evaluating the sum, and using the formula for the sum of squares, we
get

\begin{align}
m_{normal}  = \frac {\frac{K+1}{2} + \frac {K(K+1)(2K+1)}{6K}
\epsilon} {K+1+ \frac {K^2 + K}{2} \epsilon}
\end{align}

Factoring out a K+1 from the numerator and denominator, expanding in
powers of epsilon, and dropping everything after the linear term, we
get

\begin{align}
m_{normal}  &= \frac {0.5 + \frac {(2K+1)}{6} \epsilon} {1 + \frac
{K}{2}
\epsilon}\\
m_{normal}  &= (0.5 + \frac {(2K+1)}{6} \epsilon)(1 - \frac {K}{2} \epsilon)\\
m_{normal}  &= (0.5 + \epsilon(\frac {K+2}{12}))
\end{align}

To find the overall value of $m$, we take a weighted average of
$m_{normal}$ and $m_{zealots}$.

\begin{align}
m &= \alpha m_{normal} + (1-\alpha) m_{zealots}\\
m &= \alpha (0.5 + \epsilon(\frac {K+2}{12})) + (1-\alpha) (0.5)\\
m &= 0.5 + \epsilon \alpha \frac {K+2}{12}
\end{align}

Next, to check for stability, we calculate $E[\frac {dm}{dt}]$. When
this figure has the same sign as $\epsilon$, the equilibrium state
is unstable, and when it has the opposite sign, the uniform
distribution is stable.

\begin{align}
E[\frac {dm}{dt}] &= m(1-\rho_K) - (1-m)(1-\rho_0)\\
&= (0.5 + \epsilon \alpha \frac {K+2}{12}) (1 - \frac {1 + K
\epsilon}{K + 1 + \epsilon\frac {K^2 + K}{2}})\\
\notag  &- (0.5 - \epsilon \alpha \frac {K+2}{12}) (1 - \frac {1 - K
\epsilon}{K + 1 + \epsilon\frac {K^2 + K}{2}})
\end{align}

Expanding in powers of $\epsilon$ and dropping everything after the
linear terms, we get

\begin{align}
E[\frac {dm}{dt}] &= (0.5 + \epsilon \alpha \frac {K+2}{12}) (1 -
\frac {1}{K+1} - \frac {K}{2} \epsilon) \\
\notag &-(0.5 - \epsilon \alpha \frac {K+2}{12}) (1 - \frac {1}{K+1}
+ \frac {K}{2} \epsilon)\\
&= \epsilon(-\frac{K}{4} + \alpha \frac {K^2 + 2K}{12K+12}) -
\epsilon (\frac{K}{4} - \alpha \frac {K^2 + 2K}{12K+12})\\
&= \epsilon (-2\frac{K}{4} +2 \alpha \frac {K^2 + 2K}{12K + 12})
\end{align}

Note that the uniform distribution is stable whenever $\frac {K}{4}
- \alpha \frac {K^2 + 2K}{12K + 12} > 0.$  As a result, we get

\begin{align}
\alpha < \frac{\frac{K}{4}}{\frac{K^2 + 2K}{12K + 12}}\\
\alpha < \frac{3}{K + 2}
\end{align}

Figure compares this analytical result to numerical results for K =
2, 3, and 4. As can be seen, the analytics and the numerics agree.
It should be noted that the discrepancy between the highest and
lowest values of $m$ corresponding to equilibrium states is
relatively high for values of $\alpha$ even only slightly higher
than the critical value of $\frac {3}{K+2}$.

\begin{figure}
\centering
\includegraphics[width=3.0in,bb=2.0in 3.2in 6.5in 7.5in]{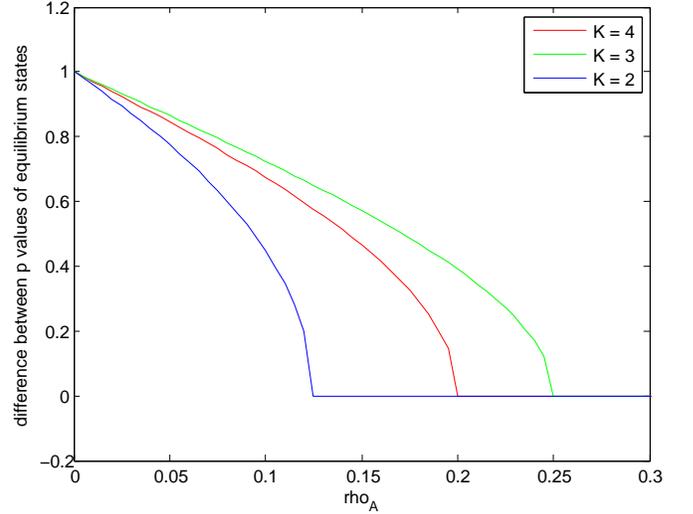}
\caption{The graph of the difference between the highest and lowest
values of $m$ for which a steady state exists in the case where
$\rho_A = \rho_B$.  For values of $\rho_A$ and $\rho_B$ slightly
below the critical value, the difference is still fairly high.}
\end{figure}

\appendices

\section{Detailed Mathematical Proofs}

In the case without zealots, $\rho_A = \rho_B = 0$.  Let the
distribution of agents in the various opinion states be represented
by the K+1 dimensional vector $[\rho_0, \rho_1, ... \rho_K]$. Let
$m$ represent the magnetization of this distribution.  Assume the
distribution to be a steady state.  Let $\phi_{i+}$ be the expected
flux over one time step from state $N_i$ to state $N_i+1$.  This is,
by definition, the expected number of listeners in state $i$ who
hear opinion A.  We define a time step such that one conversation
occurs per unit time.  The probability that a listener is in state
$i$ is simply $\rho_i$.  The probability that a listener hears
opinion A is roughly the same for all potential listeners, (assuming
a sufficiently large complete network,) and is equal to the
magnetization.  Thus, we have

\begin{align}
\phi_{i+} = m \times \rho_i
\end{align}

for $i = 0, 1, 2 ... K-1$.  $\phi_{K+} = 0$ because there is no
state $N_{K+1}$ for the agents to move into.  Similarly, we can
define $\phi_{i-}$ as the expected flux from state $N_i$ to state
$N_{i-1}$.  We have

\begin{align}
\phi_{i-} = (1-m) \times \rho_i
\end{align}

If a distribution is a steady state, the flux from state $N_i$ to
state $N_{i+1}$ must be matched by the flux from state $N_{i+1}$ to
state $N_i$.  Thus, we have the equations

\begin{align}
\phi_{(i+1)-} = \phi_{i+}\\
(1-m) \times \rho_{i+1} = m \times \rho_i\\
\frac{\rho_{i+1}}{\rho_i} = \frac {m}{1-m}
\end{align}

for $m \neq 1$.  The special cases of $m = 0$ and $m = 1$ are the
absorbing states, and are clearly steady states.  The case $m = 0.5$
produces a steady state.  In this case

\begin{align}
\frac{\rho_{i+1}}{\rho_i} = \frac {0.5}{0.5} = 1\\
\rho_{i+1} = \rho_i
\end{align}

Because all K+1 population densities are equal and they sum to 1, we
have $\rho_i = \frac {1}{K+1}$ for all i.  The only question is
whether or not the magnetization of the resulting distribution is,
in fact, 0.5.

\begin{align}
m = \frac {\sum_{i = 1}^{K} i \times \rho_i}{K}\\
m = \frac {\sum_{i = 1}^{K} i \times \frac {1}{K+1}}{K}\\
m = \frac { \frac {K \times {(K+1)}}{2} \times \frac {1}{K+1}}{K}\\
m = 0.5
\end{align}

Because the distribution is geometric and the resulting
magnetization is equal to the magnetization used to generate the
distribution, it is a steady state.  Thus, steady states are
achieved for the geometric distributions generated by magnetizations
of 0, 0.5, and 1.

In the case where $0<m<0.5$, when $m$ is used to generate a
distribution and $K \geq 2$, the resulting magnetization is strictly
less than the one used to generate the distribution, and the when
$0.5<m<1$, the resulting magnetization is strictly greater than the
one used to generate the distribution.  This means that, in either
case, the distribution is not a steady state as the ratio $\frac
{\rho_{i+1}}{\rho_i} \neq \frac {m}{1-m}$.  The proof is shown for
the case $0.5<m<1$ and an analogous proof can be used for the case
$0<m<0.5$.

Suppose K = 1.  Let $\rho_0 = 1-x$ and $\rho_1 = x$.  The
magnetization of this distribution is

\begin{align}
m = 1 \times x = x
\end{align}

and the ratio of adjacent states is

\begin{align}
\frac {\rho_1}{\rho_0} = \frac {x}{1-x} = \frac {m}{1-m}
\end{align}

Thus, for K = 1, any geometric distribution is a steady state. If K
is increased by 1 and the same ratio $\frac {x}{1-x}$ between
adjacent states and $0.5<x<1$, then the value of $m$ for the new
distribution is strictly greater than the value of $m$ for the
distribution with a lower value of K.

Consider a geometric distribution on K+1 states with $\frac
{\rho_1}{\rho_0} = \frac {x}{1-x}$.  Suppose this distribution has
magnetization $m_{old}$.  Suppose that this distribution is extended
to K+2 states by keeping the values of $\rho_0, \rho_1...\rho_K$ and
setting $\rho_{K+1} = \rho_K \times \frac {x}{1-x}$.  We calculate
the magnetization of the new distribution, $m_{new}$.

\begin{align}
m_{new} &= \frac {\sum_{i=1}^{K+1} (i \rho_i)}{\sum_{i=0}^{K+1}
\rho_i}\\
m_{new} &= \frac {\frac{K}{K+1} {m_{old}} {\sum_{i=0}^{K}
\rho_i} + \rho_{K+1}} {{\sum_{i=0}^{K+1} \rho_i}}\\
m_{new} &= m_{old} + \frac {\frac{-1}{K+1} m_{old} {\sum_{i=0}^{K}
\rho_i} + {(1-m_{old})} \rho_{K+1}} {{\sum_{i=0}^{K+1} \rho_i}}
\end{align}

We examine the sign of $\frac{-1}{K+1} \times m_{old} *
{\sum_{i=0}^{K} \rho_i} + {1-m_{old}} \times \rho_{K+1}$.  Note
that, because of the ratio $\frac {\rho_{i+1}}{\rho_i} = \frac
{x}{1-x}$ and the restrictions taken on $x$, $\rho_{K}$ must be
greater than any of the other values of $\rho_i$.  Thus, $\rho_{K} >
\frac {\sum_{i=0}^{K} \rho_i} {K+1}$, so $\rho_{K+1} > \frac {x
\times \sum_{i=0}^{K} \rho_i} {(1-x)(K+1)}$. In the case K = 1, we
have $m_{old} = x$ and we will show, by induction, that for $K \geq
2$, $m_{old} \geq x$, which gives us $\rho_{K+1} > \frac {m_{old}
\times \sum_{i=0}^{K} \rho_i} {(1-m_{old})(K+1)}$.

\begin{align}
&m_{new} = m_{old} + \frac {\frac{-1}{K+1} m_{old} {\sum_{i=0}^{K}
\rho_i} + {(1-m_{old})} \rho_{K+1}} {{\sum_{i=0}^{K+1} \rho_i}}\\
&m_{new} > m_{old} + \frac {\frac{-1}{K+1} m_{old} {\sum_{i=0}^{K}
\rho_i} + \frac {m_{old}}{1-m_{old}}{(1-m_{old})} \frac
{\sum_{i=0}^{K} \rho_i} {K+1}} {{\sum_{i=0}^{K+1} \rho_i}}\\
&m_{new}
> m_{old} + \frac {\sum_{i=0}^{K} \rho_i} {\sum_{i=0}^{K+1} \rho_i}
\times \frac {m_{old} \frac {-1}{K+1} + \frac
{1}{K+1}}{{\sum_{i=0}^{K+1} \rho_i}}\\
&m_{new} > m_{old}
\end{align}

Because $m_{new} > m_{old}$ and $m_{old} \geq \frac {x}{1-x}$ we
have $m_{new} > \frac {x}{1-x}$ so the magnetization of the new
distribution is not equal to $\frac {x}{1-x}$, so the distribution
is not a steady state.  Note that this inequality also validates the
induction used earlier.

Having found the steady states, we examine their stability. Consider
the state of consensus at B.  The vector of population densities is
$[1, 0, 0, ... 0]$.  Note that, because the sum of population
densities must equal 1, we can determine the distribution from the
vector $[\rho_1, \rho_2, ... \rho_K]$.  For consensus at B, this is
$[0, 0, ... 0]$.  Assume that the initial condition is a slight
perturbation from consensus, and $[\rho_1, \rho_2, ... \rho_K] =
[\epsilon_1, \epsilon_2, ... \epsilon_K]$.  Consider the
"magnetization norm" $||x||_m = (1/K) (|x_1| + 2|x_2| + ... +
K|x_K|)$.  First, we show that this quantity is, in fact, a norm.
Because it is a sum of absolute values, it is nonnegative, and 0
only when each component is 0.

\begin{align}
||c \times x||_m &= \frac {|c \times x_1| + 2|c \times x_2| + ... +
K|c \times x_K|}{K}\\
& = \frac {|c| \times |x_1| + 2|c| \times |x_2| + ... + K|c| \times
|x_K|}{K}\\
& = |c| ||x||_m\\
||x+y||_m &= \frac {|x_1 + y_1| + 2|x_2 + y_2| + ... + K|x_K +
y_K|}{K}\\
& \leq \frac {|x_1| + |y_1| +  ... + K|x_K| +
K|y_K|}{K}\\
& = ||x||_m + ||y||_m
\end{align}

Next, note that for realizable distributions, $x_i \geq 0$, so the
absolute value bars can be ignored.  Thus, $||x||_m = (1/K) (x_1 +
2(x_2) + ... + K(x_K)) = m$.  At consensus, $||x||_m = 0$.  We
examine what happens in a neighborhood of consensus such that
$||x||_m \le \frac {1}{2KN}$.  Note that this means that $\rho_i =
0$ for $i > K/2$ and that $\rho_0 \le 1-2m$.  We examine $\Delta m$,
the change in $m$ over a single time step.

\begin{align}
\Delta m &= \frac {m \rho_0 + (2m-1) (1-\rho_0)}{NK}\\
& \le \frac {m (1-2m) + (2m-1) (1-(1-2m))}{NK}\\
& = \frac {m - 2m^2 + (4m^2-2m)}{NK}\\
& = \frac {-m + 2m^2 + (4m^2-2m)}{NK}\\
& < 0
\end{align}

for $0<m<0.5$.  Thus, the value of $||x||_m$ will decrease, and the
point is stable.  An analogous argument holds for the point of
consensus at A.

Next, we examine the stability of the uniform distribution. In a
sufficiently small neighborhood of this distribution, none of the
population densities are 0.  This allows us to describe the
distribution with the vector of ratios of densities of adjacent
states, $[\frac {\rho_1}{\rho_0}, \frac {\rho_2}{\rho_1}, \frac
{\rho_3}{\rho_2}...\frac {\rho_K}{\rho_K-1}]$.  For the uniform
distribution, this vector is $[1, 1, 1...1]$.  Assume a small
perturbation of the form $[\epsilon, \epsilon,
\epsilon...\epsilon]$, so that the vector is $[1+\epsilon,
1+\epsilon...1+\epsilon]$, where $\epsilon > 0$.  We claim that if
the minimum ratio of densities of adjacent states is some
$r_{min}(t)>1$, then after one time step, the minimum ratio of
densities is some $r_{min}(t+1) > r_{min}(t)$.  Let's first examine
the ratio $\frac {\rho_{i+1}}{\rho_i}$ for $i = 1, 2, ...K-2$.  We
fix $\rho_i$ and attempt to minimize $\frac {\rho_{i+1}}{\rho_i}
(t+1)$.  Suppose that $\rho_{i+1} > r_{min}(t) (\rho_i)$.  Then, for
sufficiently large N, $\rho_{i+1}$ and $\rho{i}$ would change
sufficiently little to keep $\frac {\rho_{i+1}}{\rho_i} (t+1) >
r_{min}(t)$.  If $\rho_{i+1} = r_{min}(t) (\rho_i)$, we attempt to
minimize $\frac {\rho_{i+1}}{\rho_i} (t+1)$.  To do this, we
minimize $\rho_{i+1}(t+1)$ and maximize $\rho_i (t+1)$.  This is
achieved by minimizing $\rho_{i+2}(t)$ and maximizing $\rho_{i-1}
(t)$.  Obeying the constraint that $\frac {\rho_{i+1}}{\rho_i} \geq
r_{min}(t)$, we get

\begin{align}
\rho_{i-1}(t) \leq \frac {\rho_i (t)}{r_{min}(t)}\\
\rho_{i}(t) = \rho_i (t)\\
\rho_{i+1}(t) = \rho_i (t)\times {r_{min}(t)}\\
\rho_{i+2}(t) \geq \rho_i (t)\times {(r_{min}(t))}^2
\end{align}

The new values of the relevant densities are

\begin{eqnarray}
\rho_i{(t+1)} = \rho_i(t)*(\frac {N-1}{N}) + \rho_{i-1}(t)\times
m(1/N) \nonumber\\
+ \rho_{i+1}(t)\times (1-m)(1/N)\\
\rho_{i+1}{(t+1)} = \rho_{i+1}(t)*(\frac {N-1}{N}) +
\rho_{i}(t)\times
m(1/N) \nonumber\\
+ \rho_{i+2}(t)\times (1-m)(1/N)
\end{eqnarray}

The quotient of these terms is

\begin{eqnarray}
\frac {\rho_{i+1}(t+1)}{\rho_{i}(t+1)} = \nonumber\\
\frac {\rho_{i+1}(t)(\frac {N-1}{N}) + \rho_{i}(t)m(1/N) +
\rho_{i+2}(t) (1-m)(1/N)} {\rho_i(t)(\frac {N-1}{N}) +
\rho_{i-1}(t)m(1/N)
+ \rho_{i+1}(t)(1-m)(1/N)}\\
\geq \frac {\rho_i(t)(r_{min}(t)(\frac {N-1}{N})+m(1/N)+r_{min}(t)^2
 (1-m)(1/N))}{\rho_i(t) ((\frac {N-1}{N}) + \frac
{m}{N r_{min}(t)} +r_{min}(t) (1-m)(1/N))}\\
= r_{min}(t)
\end{eqnarray}

Thus, the if the minimum ratio of densities of adjacent opinion
states decreases, it must involve one of the extreme states.  We
examine the case $i = 0$.  Similar to before, we have

\begin{align}
\rho_{0}(t) = \rho_0 (t)\\
\rho_{1}(t) = \rho_0 (t)\times {r_{min}(t)}\\
\rho_{2}(t) \geq \rho_0 (t)\times {(r_{min}(t))}^2\\
\end{align}

The new values of the relevant densities are

\begin{eqnarray}
\rho_0{(t+1)} = \rho_0(t)*(\frac {N-m}{N})\nonumber\\
 + \rho_{1}(t)\times (1-m)(1/N)\\
\rho_{1}{(t+1)} = \rho_{1}(t)*(\frac {N-1}{N}) \nonumber \\
+\rho_{0}(t)\times
m(1/N) + \rho_{2}(t)\times (1-m)(1/N)\\
\end{eqnarray}

The quotient of these terms is

\begin{align}
&\frac {\rho_{1}(t+1)}{\rho_{0}(t+1)} = \frac {\rho_{1}(t)(\frac
{N-1}{N}) + \rho_{0}(t) \frac {m}{N} + \rho_{2}(t)(1-m)(1/N)}
{\rho_0(t)(\frac {N-m}{N})
+ \rho_{1}(t)(1-m)(1/N)}\\
& \geq \frac {\rho_0(t)(r_{min}(t)(\frac
{N-1}{N})+m(1/N)+r_{min}(t)^2 (1-m)(1/N))}{\rho_0(t) ((\frac
{N-m}{N}) +r_{min}(t) (1-m)(1/N))}
\end{align}

The numerator of this fraction is a special case of the numerator of
equation 58.  The denominator differs by an additive term of
$(\rho_0(t))\times (\frac {-m}{N \times r_{min}(t)} + \frac
{1-m}{N})$. Because $r_{min}(t) > 1$, we have

\begin{align}
\frac {m}{1-m} > r_{min}(t)\\
\frac {m}{r_{min}(t)} > 1-m\\
\frac {m}{N \times r_{min}(t)} > \frac {1-m}{N}\\
\frac {-m}{N \times r_{min}(t)} + \frac {1-m}{N} < 0
\end{align}

Thus, the denominator is smaller, so the quotient is larger, and we
have

\begin{align}
\frac {\rho_{1}(t+1)}{\rho_{0}(t+1)} > r_{min}(t)
\end{align}

Finally, we examine the case where $i = K-1$.

\begin{align}
\rho_{K-2}(t) \leq \frac {\rho_K-1 (t)}{r_{min}(t)}\\
\rho_{K-1}(t) = \rho_K-1 (t)\\
\rho_{K}(t) = \rho_K-1 (t)\times {r_{min}(t)}
\end{align}

The new values of the relevant densities are

\begin{eqnarray}
\rho_K-1{(t+1)} = \rho_K-1(t)*(\frac {N-1}{N}) +\nonumber\\
\rho_{K-2}(t)\times
m(1/N) + \rho_{K}(t)\times (1-m)(1/N)\\
\rho_{K}{(t+1)} = \rho_{K}(t)*(\frac {N-1+m}{N}) + \nonumber\\
\rho_{K-1}(t)\times m(1/N)
\end{eqnarray}

Comparing this quotient to equation 64, we find that the denominator
can be factored into the same form, but the numerator has an
additive constant of $\rho_{K} \times (\frac {m}{N} - \frac
{(1-m)\times (r_{min})}{N})$

Again, we have

\begin{align}
\frac {m}{1-m} > r_{min}(t)\\
m > (1-m) \times {r_{min}(t)}\\
\frac {m}{N } > \frac {(1-m) \times {r_{min}(t)}}{N}\\
\frac {m}{N } -  \frac {(1-m) \times {r_{min}(t)}}{N} > 0
\end{align}

Thus, the numerator is greater and the denominator is the same as in
equation 58, so the quotient is greater.  Again, we have
\begin{align}
\frac {\rho_{K}(t+1)}{\rho_{K-1}(t+1)} > r_{min}(t)
\end{align}

Thus, $r_{min}(t)$ in this case is a non-decreasing function of $t$,
so it will never drop back down to 1.  The uniform distribution is
unstable.

In the case of unilateral zealots, either $\rho_A = 0$ or $\rho_B =
0$.  For simplicity, we assume $\rho_A = 0$ and infer the
corresponding results for $\rho_B = 0$ by symmetry.  The
distribution of agents in each non-zealot opinion state must still
be geometric, so for fixed $\rho_B$ we can determine the
distribution by the parameter $m$.  Note that, given a distribution
of normal agents and zealots, we have the relation

\begin{align}
m = m_{normal}{(1-\rho_B - \rho_A)} + m_{zealots}{(\rho_A+\rho_B)}
\end{align}

Note that, in the case where $\rho_A = 0$, $m_{zealots} = 0$, so the
equation reduces to

\begin{align}
m = m_{normal}{(1-\rho_B)}
\end{align}

In the case without zealots, we know that the only values of $m$
which yield steady states are 0, 0.5, and 1.  The question now
arises, for other values of $m$ is there a value of $\rho_B$ such
that there is a steady state with that value of $m$.  We solve for
any value of $\rho_B$ that would make this possible.

\begin{align}
m = m_{normal}{(1-\rho_B)}\\
m = m_{normal} - m_{normal} \rho_B\\
\rho_B = \frac {m_{normal} - m}{m_{normal}}
\end{align}

We have the added restrictions that $0 \leq \rho_B \leq 1$.  This
can be achieved whenever ${m_{normal} - m} \geq 0$.  As we have
previously shown, this occurs whenever $0 < m < 0.5$.  Thus, for
those values of $m$ there exists a unique $\rho_B$ such that
equilibrium exists at that particular $m$.  This is consistent with
what has been shown in the case K = 2.  In that particular case, it
was found that for sufficiently small $\rho_B$ there were three
values of $m$ where equilibrium was achieved.  One was consensus at
B.  The others started at $m = 0.5$ and $m = 1$ and as $\rho_B$
increased, the corresponding equilibrium values of $m$ approached
each other.  At some critical value of $\rho_B$, the equilibrium
states meet, and above that critical value, they vanish, and only
the one equilibrium point, (consensus at B,) remains.  This will
occur if $\frac {d^2 \rho_B}{dm^2} < 0$ over the interval $(0.5,
1)$.

The state of consensus at B can be shown to be stable using the
aforementioned magnetization norm.  Again, examine the neighborhood
of consensus such that $||x||_m < \frac {1}{2KN}$.

\begin{align}
\Delta m &= \frac {m (\rho_0 + \rho_B) + (2m-1) (1-\rho_0 - \rho_B)}{NK}\\
& \le \frac {m (1-2m) + (2m-1) (1-(1-2m))}{NK}\\
& = \frac {m - 2m^2 + (4m^2-2m)}{NK}\\
& = \frac {-m + 2m^2 + (4m^2-2m)}{NK}\\
& < 0
\end{align}

over the interval $(0, 0.5)$, showing that the point is stable.
assuming that the pattern of three equilibrium points for
sufficiently small $\rho_B$ holds, it can be shown that the middle
one is unstable.  Note that, under these assumptions, the value of
$\rho_B$ needed to maintain equilibrium increases as $m$ increases.
Thus, as $m$ increases and is used to generate a geometric
distribution, the resulting value of $m$ will be greater than the
one used to generate the distribution.  (It should be noted that in
the case of the last equilibrium point, the opposite is true.)
Because the quantities $\rho_0, \rho_1...$ are nonzero, and because
the value of $\rho_B$ is fixed, we can determine the distribution
(and any others near it,) with the vector $[r_1, r_2, ... r_K]$,
where $r_i = \frac {\rho_i}{\rho_{i-1}}$.  At equilibrium, there is
a magnetization value $m_{eq}$, and the vector of density ratios is
$[\frac {m}{1-m}, \frac {m}{1-m}, \frac {m}{1-m}, ... \frac
{m}{1-m}]$.  We assume that the initial state is perturbed slightly
from the equilibrium state, and the vector of density ratios is now
$[\frac {m^*}{1-m^*}, \frac {m^*}{1-m^*}, ... \frac {m^*}{1-m^*}]$,
where $m^* = m_{eq} + \epsilon$.  We assume that if the minimum
value of the ratio of densities of adjacent opinion states is no
less than $\frac {m^*}{1-m^*}$, then after one time step the minimum
ratio will still be no less than $\frac {m^*}{1-m^*}$, near enough
the equilibrium point.  Consider the ratio between two ambivalent
states.  Just like the case without zealots, we get

\begin{align}
\rho_i{(t+1)} = \rho_i(t)*(\frac {N-1}{N}) + \rho_{i-1}(t)\times
m(1/N)\nonumber\\
 + \rho_{i+1}(t)\times (1-m)(1/N)\\
\rho_{i+1}{(t+1)} = \rho_{i+1}(t)*(\frac {N-1}{N}) +
\rho_{i}(t)\times m(1/N)\nonumber\\
 + \rho_{i+2}(t)\times (1-m)(1/N)
\end{align}

And the quotient is still bounded from below by $\frac
{m^*}{1-m^*}$.

Let's consider the ratio between $\rho_1$ and $\rho_0$.  Because of
the assumption on the distribution, $\rho_1 \ge \frac {m^*
\rho_0}{1-m^*}$.  As before, the only way that the ratio could be
smaller than $\frac {m^*}{1-m^*}$ after one time step for
sufficiently large N is if $\rho_1 = \frac {m^* \rho_0}{1-m^*}$.
Furthermore, recall that $\rho_2 \ge \frac {m^* \rho_1}{1-m^*}$.
The ratio of $\rho_1$ to $\rho_0$ after one time step is

\begin{align}
\frac {\rho_{1}(t+1)}{\rho_{0}(t+1)} = \frac {\rho_{1}(t)(\frac
{N-1}{N}) + \rho_{0}(t)(m/N) + \rho_{2}(t)(1-m)(1/N)}
{\rho_0(t)(\frac {N-m}{N}) + \rho_{1}(t)(1-m)(1/N)}\\
\geq
\frac {\rho_0(t)(\frac {m^*}{1-m^*}*(\frac {N-1}{N})+m(1/N)+\frac
{m^*}{1-m^*}^2 \times (1-m)(1/N))}{\rho_0(t) ((\frac {N-m}{N})
+\frac {m*}{1-m*} \times (1-m)(1/N))}
\end{align}

Note that, for the greater value of $m$, a greater value of $\rho_B$
is needed to bring equilibrium.  Thus, if $m^*$ is above the
equilibrium value, than the value of $m$ of the distribution
generated by $m^*$ is greater than $m^*$.  With this, and logic
similar to the no-zealot case, it can be shown that the ratio of
$\rho_1$ to $\rho_0$ will remain above $\frac {m^*}{1-m^*}$.  An
analogous argument holds for the ratio of $\rho_K$ to $\rho_{K-1}$.
Thus, the point is unstable.

\section{acknowledgements}
\textbf{Acknowledgements:} This work was
supported in part by the Army Research Office Grants No.
W911NF-09-1-0254 and W911NF-12-1-0546. The views and conclusions
contained in this document are those of the authors and should not
be interpreted as representing the official policies, either
expressed or implied, of the Army Research Office or the U.S.
Government.

%


\newpage
\bibliographystyle{IEEEtran}
\bibliography{IEEEabrv,project}

[1] J. Xie, S. Sreenivasan, G. Korniss, W. Zhang, C. Lim, and B. K.
Szymanski, 2001, Social consensus through the influence of committed
minorities, Phys. Rev. E 84, 011130\newline

[2] M Starnini, A. Baronchelli, and R. Pastor-Sartorras, 2012,
Ordering dynamics of the multi-state voter model, Journal-ref: J.
Stat. Mech, P. 10027.\newline

[3] A. Baronchelli, V. Loreto, and L. Steels, 2008, In-depth
analysis of the Naming Game dynamics: the homogeneous mixing case,
Int. J. Mod. Phys. C 19, 785\newline

[4] W. Zhang, C. Lim, S. Sreenivasan, J. Xie, B. K. Szymanski, and
G. Korniss, 2011, Social influencing and associated random walk
models: Asymptotic consensus times on the complete graph, Chaos 21,
025115\newline

[5] X. Castelló, A. Baronchelli, and V. Loreto, 2009, Consensus and
ordering in language dynamics, European Physical Journal B 71,
557-564\newline

[6] A. Baronchelli, M. Felici, V. Loretto, E. Caglioti, and L.
Steels, 2005, Sharp transition towards shared vocabularies in
multi-agent systems, J. Stat Mech:Theory Exp P06014: 0509075\newline

[7] A. Baronchelli, L Dall'Asta, A Barrat, and V. Loreto, 2007
Nonequilibrium phase transition in negotiation dynamics. Phys Rev E
76: 051102\newline

[8] Q. Lu, G Korniss, and B Szymanski, 2009, The naming game in
social networks: community formation and consensus engineering. J
Econ Interact Coord 4: 221\newline

[9] M. Mobilia, 2003, Does a single zealot affect an infinite group
of voters? Phys Rev Lett 91: 028701\newline

[10] S. Galam and F. Jacobs, 2007, The role of inflexible minorities
in the breaking of democratic opinion dynamics. Physica A 381: 366 -
376\newline

[11] M. Mobilia, A. Petersen, and S. Redner, 2007, On the role of
zealotry in the voter model. J Stat Mech:Theory Exp 2007:
P08029.\newline

[12] E. Yildiz, D. Acemoglu, A. Ozdaglar, A. Saberi, and A.
Scaglione, 2011, Discrete Opinion Dynamics with Stubborn Agents.
Operations Research manuscript OPRE-2011-01-026\newline

[13] W. Zhang, C. Lim, and B. Szymanski, Analytic Treatment of
Tipping Points for Social Consensus in Large Random Networks, 2012,
Physical Review E, 86(6) 061134\newline

[14] N. Lanchier, The Naming Game in Language Dynamics Revisited,
2013, arXiv:1301.0124\newline

[15] B. Lin, J. Ren, H. Yang, and B. Wang, Naming Game on
small-world networks: the role of clustering structure, 2006,
arXiv:physics/0607001\newline

[16] B. Bhattacherjee, S. Manna, and A. Mukherjee, Information
sharing and sorting in a community, 2013, arXiv:1306.0646\newline

[17] J. Xie, J. Emenheiser, M. Kirby, S. Sreenivasan, B. K.
Szymanski, and G. Korniss, Evolution of opinions on social networks
in the presence of competing committed groups, 2012, arXiv:1112.6414
[physics.soc-ph]

[18] Y. Treitman, C. Lim, W. Zhang, and A. Thompson, Naming Game
with Greater Stubbornness and Unilateral Zealots, 2013, IEEE
$2^{nd}$ annual International Network Science Workshop, 2013, print.

[19] G. Verma, A. Swami, K. Chan, The Effect of Zealotry in the
Naming Game Model of Opinion Dynamics, MILITARY COMMUNICATIONS
CONFERENCE, 2012 - MILCOM 2012 Orlando, Oct. 29 2012-Nov. 1 2012,
IEEE Computer Science Press

[20] D. J. D. Earn, and S. A. Levin, Global asymptotic coherence in
discrete dynamical systems. Proceedings of the National Academy of
Sciences, USA 103 (11): 3968-3971,  2006

[21] A. P. Kinzig , P. R. Ehrlich , L. J. Alston , K. Arrow , S.
Barrett , T. G. Buchman , G. C. Daily , B. Levin , S. Levin , M.
Oppenheimer , E. Ostrom and D. Saari, Social Norms and Global
Environmental Challenges: The Complex Interaction of Behaviors,
Values, and Policy, BioScience, 63(3):164-175. 2013.
http://www.bioone.org/doi/full/10.1525/bio.2013.63.3.5

[22] M. Turalska1, B. J. West, and P. Grigolini,  Role of committed
minorities in times of crisis, SCIENTIFIC REPORTS | 3 : 1371 | DOI:
10.1038/srep01371

\end{document}